\documentclass [12pt, draft]{article}
\usepackage{amsmath, amsfonts, amsthm, amssymb}
\newtheorem{theorem}{Theorem}[section]
\newtheorem{prove}{Proof of Theorem}[section]
\newtheorem{lemma}[theorem]{Lemma}
\newtheorem{definition}{Definition}[subsection]
\newtheorem{proposition}[definition]{Proposition}

\begin{document}
\title{\bf Classification of 5-dimensional MD-algebras having commutative derived ideals}
\author{${\bf Le\,\, Anh\,\, Vu}^{*}$ \, and \, ${\bf Kar\,\, Ping\,\, Shum}^{**}$
\\${}^{*}$ \footnotesize{Department of Mathematics and Informatics}
\\\footnotesize{University of Pedagogy, Ho Chi Minh City, Viet nam}
\\\footnotesize{E-mail:\, vula@math.hcmup.edu.vn}
\\${}^{**}$ \footnotesize{Faculty of Science, The Chinese University of Hong Kong}
\\\footnotesize{E-mail:\, kpshum@maths.hku.hk}}
\date{}
\maketitle

\begin{abstract}
    In this paper, we study  a subclass  of  the class of  MD-algebras,
    i.e., the class of solvable real Lie algebras such that the K-orbits of
    its corresponding connected and simply connected Lie groups are either
    orbits of dimension zero or orbits with maximal dimensions. Our main result is
    to classify, up to isomorphism, all the 5-dimensional MD-algebras having commutative
    derived ideals.
\end{abstract}



\subsection*{Introduction}
The concept of C*-algebras was first introduced by Gelfand and
Naimark in 1943. It is well known that C*-algebras can be applied to
mathematics, mechanics and physics, however, the problem of
describing the structure of C*-algebras, in general, is still
open.\\

The method of describing the structure of C*-algebras by using
K-functors was first suggested by D. N. Diep ([2]) in 1974.  By
applying the K-homology functors proposed by  Brown - Douglas -
Fillmore (for brevity, the BDF K-functors), Diep gave a description
for the C*(Aff$\mathbb{R}$) of the group Aff$\mathbb{R}$ of the
affine transformations of the real line. In 1975, by using the
method of Diep, J. Rosenberg ([7], [8]) gave a description for the
C*-algebra of the group Aff$\mathbb{C}$ and some other groups. In
1977, D.N.Diep ([3]) further gave a complete system of invariants of
C*-algebras of type I by using the  BDF K-homology functors. Hence,
it is natural to propose the following two general problems: \\
\begin{itemize}
    \item[$\bullet$] Generalize the K-homology functors so that these functors
    can be applied to a larger class of C*-algebras.\\
    \item[$\bullet$] Find the C*-algebras which can be described by
    using the generalized K-functors.\\
\end{itemize}

Concerning the first problem, we note that G. G. Kasparov ([5]) in
1980 introduced the concept of  KK-functors which is a generalized
concept of BDF K-homology functors. Then by using  KK-functors, G.G.
Kasparov described the C*-algebra of the Heisenberg groups
$H_{2n+1}$.\\

For the second problem, it was noticed that this problem is closely
related with the Orbit Method proposed by A.A. Kirillov ([6]) in
1962. After studying the Kirillov's Orbit Method, Diep in 1980
suggested to consider the class of Lie groups and Lie algebras MD
and $\overline{MD}$ ([4]) so that the C*-algebras of them can be
described by using  KK-functors. If $G$ is an n-dimensional real Lie
group, then $G$ is called a MDn-group or a MD-group with dimension n
iff the orbits of $G$ in the K-representation (K-orbits) are orbits
of dimension zero or orbits of maximal dimension (i.e. dimension k,
where k is some even constant, $k \leq n$). When k = n, we call $G$
an $\overline{MDn}$-group or $\overline{MD}$-group of dimension n.
The corresponding Lie algebra Lie(G) of G is said to be an
MDn-algebra or $\overline{MDn}$-algebra, respectively. It is clear
that the class $\overline{MD}$ is a subclass of the class MD. Thus,
the problem of classifying MD-algebras, describing the
K-representation of MD-groups and characterizing the C*-algebras of
MD-groups is significant. Note that all the Lie algebras and the Lie
groups of dimension $n$ with $n < 4$ are MD-algebras and MD-groups,
and moreover they can be listed easily. So we only take interest in
MDn-groups and MDn-algebras for $n \geq 4$.\\

We remark here that all $\overline{MD}$-algebras (of arbitrary
dimension) was classified, up to  isomorphism, by H. H. Viet in [9].
This class includes only the following algebras:\\
\begin{itemize}
    \item[$\bullet$]$\mathbb{R}^n$ - The commutative Lie Algebra of dimension
    n;\\
    \item[$\bullet$]Lie(Aff$\mathbb{R}$) - The Lie algebra of the group
    of affine traformations of the real straight line;\\
    \item[$\bullet$]Lie(Aff$\mathbb{C}$) - The Lie algebra of the group of affine
transformations of the complex straight line.\\
\end{itemize}
It is noteworthy that Viet [9] also described the C*-algebras of the
universal covering of group Aff$\mathbb{C}$ by using KK-functors.
Thus, the C*-algebras of all groups of the class $\overline{MD}$
were
described by Diep, Rosenberg and  Viet.\\

The problem for the class of MD-algebras is much more complicated
than $\overline{MD}$-algebras. In 1984, Dao Van Tra [11] listed all
MD4-algebras. In 1990, all MD4-algebras were classified, up to
isomorphism, by Vu (see [12], [13], [14]). Until quite recently, Vu
together with Nguyen Cong Tri, Duong Minh Thanh and Duong Quang Hoa
introduced some MD5 - algebras and MD5 - groups (see [15], [16],
[17], [18], [19], [20]). Until the present moment, there is no
complete classification for MDn-algebras with $n \geq 5$.\\

On the other hand, by studying the foliated manifold, Connes ([1])
in 1982 proposed the notion of C*-algebras associated with a
measured foliation. The following question naturally arises: Can we
describe the Connes C*-algebras by using KK-functors? In fact, Torpe
has shown in [10] that the KK-functors are very useful and effective
to describe the structure of Connes C*-algebras associated with the
Reeb foliations.\\

The other reason for studying the class MD is based on the following
fact: if G is a certain MD-group, then the family of its K-orbits
with maximal dimension forms a measured foliation. This foliation is
called MD-foliation associated with G. Furthermore, the C*-algebra
of G can be easily described when the Connes C*-algebra of
MD-foliation associated with G is known. Hence, the problem of
classifying the topology and describing the Connes
C*-algebras of the class of MD-foliations is worth to study.\\

On this aspect, Vu in 1992 gave a topological classification of all
MD4-foliations and described all Connes C*-algebras of them by using
the KK-functors (see [12], [13], [14]). We noticed that the Connes
C*-algebras of MDn-foliations with $n > 4$ has not yet been
described.  Following [9],  if $\mathcal{G}$ is an MD-algebra then
the second derived ideal $\mathcal{G}^2 = [\mathcal{G}^1,
\mathcal{G}^1] = [[\mathcal{G}, \mathcal{G}],[\mathcal{G},
\mathcal{G}]]$ is commutative, however, the converse is not true.
Therefore, we need to consider only $\mathcal{G}$ for which
$\mathcal{G}^2$ is commutative. In particular, if $\mathcal{G}^2$ =
0 (i.e. $\mathcal{G}^1$ is commutative) then $\mathcal{G}$ could be
an MD-algebra. Hence, we will restrict ourself only to this case.
Our main result is to classify, up to an isomorphism, all
MD5-algebras $\mathcal{G}$ having commutative derived ideal
$\mathcal{G}^1 = [\mathcal{G}, \mathcal{G} ]$. The topology of
MD5-foliations associated with the MD5-groups and the description of
Connes C*-algebras of these foliations will be considered and
studied later on.\\

\section{Preliminaries}
We first recall in this Section some preliminary results and
notations which will be used in the sequel. For more detailed
information, the reader is referred to [4] and [6].

\subsection{The co-adjoint Representation and  K-orbits}
    Let G be a Lie group. Let $\mathcal{G}$ = Lie(G) be the Lie algebra of G and we use
$\mathcal{G^{*}}$ to denote the dual space of $\mathcal{G}$. For every g $\in$
G, we denote the internal automorphism associated with g by
$A_{(g)}$, and whence,  $ A_{(g)}: G \longrightarrow G $ can be defined as
follows
$$A_{(g)}(x):=\, g.x.g^{-1},\, \forall x \in G.$$\\
The above automorphism induces the following mapping:
$${A_{(g)}}_{*}:\mathcal{G}\longrightarrow \mathcal{G}\qquad \qquad \qquad \qquad $$
$$\qquad\qquad \qquad \qquad\qquad X\longmapsto {A_{(g)}}_{*}(X):\, =\;\frac{d}{dt}
[g.exp(tX)g^{-1}]\mid_{t=0}$$\\ which is called \textit {the tangent
mapping} of $A_{(g)}$.\\

\ We now formulate the following definitions.\\

\begin{definition}The action
$$Ad :G\longrightarrow Aut(\mathcal{G})$$
$$\qquad\qquad\qquad g\longmapsto Ad(g):\, =\;  {A_{(g)}}_{*}$$
is called the adjoint representation of G in $\mathcal{G}$.\\
\end{definition}

\begin{definition}The action
$$K:G\longrightarrow Aut(\mathcal{G}^{*})$$
$$g\longmapsto K_{(g)}$$
such that
$$\langle K_{(g)}F,X\rangle :\,=\langle F, Ad(g^{-1})X\rangle ;
\quad(F\in {\mathcal{G}}^{*},\, X \in \mathcal{G})$$ is called the
co-adjoint representation or K-representation of G in $\mathcal{{G}^{*}}$.\\
\end{definition}

\begin{definition} Each orbit of the co-adjoint representation
of G is called a K-orbit of G.\\
\end{definition}

Thus, for every $F \in \mathcal{G}^{*}$, the K-orbit containing $F$
defined above can be written by
$${\varOmega}_{F}:= \{K_{(g)}F / g \in G \}.$$

The dimension of every K-orbit of an arbitrary Lie group G is always
even. In order to define the dimension of the K-orbits
${\varOmega}_{F}$ for each F from the dual space $\mathcal{G^{*}}$
of the Lie algebra $\mathcal{G}$ = Lie (G) of G, it is useful to
consider the following skew-symmetric bilinear form $B_{F}$ on
$\mathcal{G}$
$$ B_{F}(X, Y) := \langle F, [X, Y]\rangle; \, \forall\, X, Y \in
\mathcal{G}.$$

Denote the stabilizer of $F$ under the co-adjoint representation of
G in $\mathcal{{G}^{*}}$ by $G_{F}$ and ${\mathcal{G}}_{F}$ :=
Lie($G_{F}$).\\

We shall need in the sequel the following result.\\
\begin{proposition} [see {[6, Section 15.1]}] $KerB_{F} = {\mathcal{G}}_{F}$
and $dim{\varOmega}_{F} = dim\mathcal{G} - dim{\mathcal{G}}_{F}.$
\hfill{$\square$}\\
\end{proposition}

\subsection{MDn-Groups \,and\,  MDn-Algebras}
\begin{definition}[{see [4, Chapter 4, definition 1.1]}]An MDn-group
is an n-dimensional real solvable Lie group such that its K-orbits
are orbits of dimension zero or maximal dimension. The Lie algebra
of an MDn-group is called an MDn-algebra. \\
\end{definition}

The following proposition gives a  necessary condition for a Lie
algebra belonging to the class of all MD-algebras.\\

\begin{proposition}[{see [9, Theorem 4]}] Let $\mathcal{G}$ be an MD-algebra.
Then its second derived ideal ${\mathcal{G}}^{2} := [[\mathcal{G},
\mathcal{G}], [\mathcal{G}, \mathcal{G}]]$ is commutative.
\hfill{$\square$}\\
\end{proposition}

We point out here that the converse of the above result is in general not true.
In other words, the above necessary condition is not a sufficient
condition. We now only consider the 5-dimensional Lie algebras
$\mathcal{G}$ having a second derived ideal $\mathcal{G}^{2} =
\{0\}$, i.e., the derived ideal $\mathcal{G}^{1}$ is commutative.
Thus, the $\mathcal{G}$ could be an MD5-algebra.\\

\section{The Main Result}
From now on, we use $\mathcal{G}$ to denote an Lie algebra of
dimension 5. We always choose a suitable basis $(\nobreak X_{1},
X_{2}, X_{3}, X_{4}, X_{5} \nobreak)$ in $\mathcal{G}$ so that
$\mathcal{G}$ is isomorphic to ${\mathbb{R}}^{5}$ as a real vector
space. The notation ${\mathcal{G}}^{*}$ will  be used to denote the
dual space of $\mathcal{G}$. Clearly,  ${\mathcal{G}}^{*}$ can be
identified with ${\mathbb{R}}^{5}$ by fixing in it the basis
$(\nobreak X_{1}^{*}, X_{2}^{*}, X_{3}^{*}, X_{4}^{*},
X_{5}^{*}\nobreak )$ which is the dual of the basis $(\nobreak
X_{1}, X_{2}, X_{3}, X_{4}, X_{5} \nobreak)$.\\

\begin{theorem}
Let $\mathcal{G}$ be an MD5-algebra whose $\mathcal{G}^{1}:$=
$[\mathcal{G}, \mathcal{G}]$ is commutative. Then the following
assertions hold.\\
\begin{itemize}
\item[I.] If $\mathcal{G}$ is decomposable, then $\mathcal{G}
\cong \mathcal{H} \oplus {\mathbb{R}}$, where $\mathcal{H}$ is an
MD4-algebra.\\

\item[II.] If $\mathcal{G}$ is indecomposable, then
we can choose a suitable basis $( X_{1}, X_{2}, X_{3},$ $X_{4},
X_{5} )$ of \, $\mathcal{G}$ such that $\mathcal{G}$ is isomorphic
to one and only one of the following Lie algebra.

    \begin{itemize}
    \item[1.]${\mathcal{G}}^{1} = \mathbb{R}.X_{5}\equiv
    \mathbb{R}.$\\
        \begin{itemize}
        \item[] $\mathcal{G}_{5,1}: [X_{1}, X_{2}]=[X_{3},X_{4}]=X_{5}$; the
others Lie Brackets are trivial.\\
        \end{itemize}
    \item[2.]${\mathcal{G}}^{1} = \mathbb{R}.X_{4} \oplus \mathbb{R}.X_{5}\equiv
{\mathbb{R}}^{2}$\\
        \begin{itemize}
        \item[2.1.] $\mathcal{G}_{5,2,1}: [X_{1}, X_{2}]=X_{4},\,\,[X_{2}, X_{3}]=X_{5}$; the others
Lie brackets are trivial.\\
        \item[2.2.] $\mathcal{G}_{5,2,2(\lambda)}: [X_{1}, X_{2}]=[X_{3},X_{4}]=X_{5},\,\, [X_{2}, X_{3}]=\lambda X_{4},\\ \lambda \in
\mathbb{R}\backslash \{0\}$; the others Lie Brackets are
trivial.\\
        \end{itemize}

    \item[3.]${\mathcal{G}}^{1} = \mathbb{R}.X_{3} \oplus \mathbb{R}.X_{4} \oplus \mathbb{R}.X_{5}
\equiv {\mathbb{R}}^{3}$,\,$ad_{X_{1}} = 0 $,\, $ad_{X_{2}} \in
End({\mathcal{G}}^{1}) \equiv Mat_{3}(\mathbb{R})$;\,
$[X_{1},X_{2}]= X_{3}$.\\

        \begin{itemize}
        \item[3.1.]${\mathcal{G}}_{5,3,1({\lambda}_{1}, {\lambda}_{2})}:$
$$ ad_{{X}_2} = \begin{pmatrix} {{\lambda}_1}&0&0\\
0&{{\lambda}_2}&0\\0&0&1 \end{pmatrix}; \quad {\lambda}_1,
{\lambda}_2 \in \mathbb{R}\setminus \lbrace 1\rbrace, \, {\lambda}_1
\neq {\lambda}_2 \neq 0 .$$ \\

        \item[3.2.]${\mathcal{G}}_{5,3,2(\lambda)}:$
$$ad_{{X}_2} = \begin{pmatrix} 1&0&0\\
0&1&0\\0&0&{\lambda} \end{pmatrix}; \quad {\lambda} \in
\mathbb{R}\setminus \lbrace 0, 1 \rbrace .$$ \\

        \item[3.3.]${\mathcal{G}}_{5,3,3(\lambda)}:$
$$ad_{{X}_2} = \begin{pmatrix} {\lambda}&0&0\\ 0&1&0\\0&0&1 \\ \end{pmatrix};
\quad {\lambda} \in \mathbb{R}\setminus \lbrace 1 \rbrace . $$\\

        \item[3.4.]${\mathcal{G}}_{5,3,4}:$
$$ad_{{X}_2} = \begin{pmatrix} 1&0&0\\
0&1&0\\0&0&1 \end{pmatrix}.$$\\

        \item[3.5.]${\mathcal{G}}_{5,3,5(\lambda)}:$
 $$ad_{{X}_2} = \begin{pmatrix} {\lambda}&0&0\\
 0&1&1\\0&0&1 \end{pmatrix}; \quad{\lambda} \in \mathbb{R}\setminus \lbrace 1 \rbrace .$$ \\

        \item[3.6.]${\mathcal{G}}_{5,3,6(\lambda)}:$
 $$ad_{{X}_2} = \begin{pmatrix} 1&1&0\\
 0&1&0\\0&0&{\lambda} \end{pmatrix}; \quad
 {\lambda} \in \mathbb{R}\setminus \lbrace 0, 1 \rbrace .$$ \\

        \item[3.7.]${\mathcal{G}}_{5,3,7}:$
 $$ad_{{X}_2} = \begin{pmatrix} 1&1&0\\
 0&1&1\\0&0&1 \end{pmatrix}.$$ \\

        \item[3.8.]${\mathcal{G}}_{5,3,8(\lambda, \varphi)}:$
 $$ad_{{X}_2} = \begin{pmatrix} cos{\varphi}&-sin{\varphi}&0\\
 sin{\varphi}&cos{\varphi}&0\\0&0&\lambda \end{pmatrix}; \quad
 \lambda \in \mathbb{R}\setminus \lbrace 0\rbrace, \, \varphi \in (0, \pi) .$$\\
        \end{itemize}

    \item[4.]${\mathcal{G}}^{1} = \mathbb{R}.X_{3} \oplus \mathbb{R}.X_{3} \oplus
\mathbb{R}.X_{4} \oplus \mathbb{R}.X_{5} \equiv {\mathbb{R}}^{4}$,\,
$$ad_{X_{1}} \in End({\mathcal{G}}^{1}) \equiv
Mat_{4}(\mathbb{R}).$$\\

        \begin{itemize}
        \item[4.1.]${\mathcal{G}}_{5,4,1({\lambda}_{1}, {\lambda}_{2}, {\lambda}_{3})}:$
 $$ad_{{X}_1} = \begin{pmatrix} {{\lambda}_1}&0&0&0\\
 0&{{\lambda}_2}&0&0\\0&0&{\lambda}_{3}&0\\0&0&0&1\end{pmatrix};$$
  ${\lambda}_1, {\lambda}_2, {\lambda}_3 \in \mathbb{R}\setminus
  \lbrace 0, 1\rbrace,\quad {\lambda}_1 \neq {\lambda}_2 \neq {\lambda}_3 \neq
  {\lambda}_1.$\\

        \item[4.2.]${\mathcal{G}}_{5,4,2({\lambda}_{1}, {\lambda}_{2})}:$
$$ad_{{X}_1} = \begin{pmatrix} {{\lambda}_1}&0&0&0\\
0&{{\lambda}_2}&0&0\\0&0&1&0\\0&0&0&1\end{pmatrix};\quad
{\lambda}_{1}, {\lambda}_{2} \in \mathbb{R}\setminus \lbrace 0, 1
\rbrace , {\lambda}_1 \neq {\lambda}_2. $$ \\

        \item[4.3.]${\mathcal{G}}_{5,4,3(\lambda)}:$
 $$ad_{{X}_1} = \begin{pmatrix}
 {\lambda}&0&0&0\\0&{\lambda}&0&0\\0&0&1&0\\0&0&0&1 \end{pmatrix}; \quad
 {\lambda} \in \mathbb{R}\setminus \lbrace 0, 1 \rbrace .$$

        \item[4.4.]${\mathcal{G}}_{5,4,4(\lambda)}:$
$$ad_{{X}_1} = \begin{pmatrix} {\lambda}&0&0&0\\0&1&0&0\\
0&0&1&0\\0&0&0&1 \end{pmatrix};\quad {\lambda} \in
\mathbb{R}\setminus \lbrace 0, 1 \rbrace.$$ \\

        \item[4.5.]${\mathcal{G}}_{5,4,5}:$
$$ad_{{X}_1} = \begin{pmatrix} 1&0&0&0\\0&1&0&0\\
0&0&1&0\\0&0&0&1 \end{pmatrix}.$$ \\

        \item[4.6.]${\mathcal{G}}_{5,4,6({\lambda}_{1}, {\lambda}_{2})}$ :
$$ad_{{X}_1} = \begin{pmatrix} {{\lambda}_1}&0&0&0\\
0&{{\lambda}_2}&0&0\\0&0&1&1\\0&0&0&1\end{pmatrix};\quad
{\lambda}_{1}, {\lambda}_{2} \in \mathbb{R}\setminus \lbrace 0, 1
\rbrace , {\lambda}_1 \neq {\lambda}_2.$$\\

        \item[4.7.]${\mathcal{G}}_{5,4,7(\lambda)}:$
$$ad_{{X}_1} = \begin{pmatrix}
{\lambda}&0&0&0\\0&{\lambda}&0&0\\0&0&1&1\\0&0&0&1 \end{pmatrix};
\quad {\lambda} \in \mathbb{R}\setminus \lbrace 0, 1 \rbrace .$$

        \item[4.8.]${\mathcal{G}}_{5,4,8(\lambda)}:$
$$ad_{{X}_1} = \begin{pmatrix}
{\lambda}&1&0&0\\0&{\lambda}&0&0\\0&0&1&1\\0&0&0&1 \end{pmatrix};
\quad {\lambda} \in \mathbb{R}\setminus \lbrace 0, 1 \rbrace .$$

        \item[4.9.]${\mathcal{G}}_{5,4,9(\lambda)}:$
$$ad_{{X}_1} = \begin{pmatrix}
{\lambda}&0&0&0\\0&1&1&0\\0&0&1&1\\0&0&0&1 \end{pmatrix}; \quad
{\lambda} \in \mathbb{R}\setminus \lbrace 0, 1\rbrace .$$\\

        \item[4.10.]${\mathcal{G}}_{5,4,10}:$
$$ad_{{X}_1} = \begin{pmatrix} 1&1&0&0\\0&1&1&0\\
0&0&1&1\\0&0&0&1 \end{pmatrix}.$$\\

        \item[4.11.]${\mathcal{G}}_{5,4,11({\lambda}_{1}, {\lambda}_{2},\varphi)}:$
$$ad_{{X}_1} = \begin{pmatrix} cos\varphi&-sin\varphi&0&0\\
sin\varphi&cos\varphi&0&0\\0&0&{\lambda}_{1}&0\\0&0&0&{\lambda}_{2}\end{pmatrix};$$
$${\lambda}_{1}, {\lambda}_{2} \in \mathbb{R}\setminus \lbrace 0 \rbrace ,
{\lambda}_1 \neq {\lambda}_2,\varphi \in (0,\pi).$$\\

        \item[4.12.]${\mathcal{G}}_{5,4,12(\lambda, \varphi)}:$
$$ad_{{X}_1} = \begin{pmatrix} cos\varphi&-sin\varphi&0&0\\
sin\varphi&cos\varphi&0&0\\0&0&\lambda&0\\0&0&0&\lambda\end{pmatrix};\quad
\lambda \in \mathbb{R}\setminus \lbrace 0 \rbrace,\varphi \in (0,\pi).$$ \\

        \item[4.13.]${\mathcal{G}}_{5,4,13(\lambda, \varphi)}:$
$$ad_{{X}_1} = \begin{pmatrix} cos\varphi&-sin\varphi&0&0\\
sin\varphi&cos\varphi&0&0\\0&0&\lambda&1\\0&0&0&\lambda\end{pmatrix};\quad
\lambda \in \mathbb{R}\setminus \lbrace 0 \rbrace,\varphi \in (0,\pi).$$ \\

        \item[4.14.]${\mathcal{G}}_{5,4,14(\lambda, \mu, \varphi)}:$
$$ad_{{X}_1} = \begin{pmatrix} cos\varphi&-sin\varphi&0&0\\
sin\varphi&cos\varphi&0&0\\0&0&\lambda&-\mu\\0&0&\mu&\lambda\end{pmatrix};$$
$$\lambda, \mu \in \mathbb{R}, \mu > 0, \varphi \in (0,\pi).$$\\
        \end{itemize}
    \end{itemize}
\end{itemize}
\end{theorem}

In proving Theorem 2.1,  we need some lemmas.\\

\begin{lemma} For $X, Y \in \mathcal{G} \backslash {\mathcal{G}}^{1},\, X \neq
Y$, by considering $ad_{X}, ad_{Y}$ as operators on
$\mathcal{G}^{1}$ we have $ad_{X}\circ ad_{Y} = ad_{Y}\circ ad_{X}$.
\end{lemma}
\begin{proof}  By using the Jacobi
identity for X, Y and consider an arbitrary element $Z \in
\mathcal{G}^{1}$, we have
\begin{align} \notag &\quad [[X, Y], Z]
+ [[Y, Z], X] + [[Z, X], Y] = 0\\
\notag \Leftrightarrow & \quad [X, [Y, Z]] - [Y, [X, Z]] = 0\\
\notag \Leftrightarrow & \quad ad_{X}\circ ad_{Y} (Z) =
ad_{Y}\circ ad_{X} (Z);\, \forall Z \in \mathcal{G}^{1}\\
\notag \Leftrightarrow & \quad ad_{X}\circ ad_{Y} = ad_{Y}\circ
ad_{X}. \end{align} \end{proof}

\vskip0.8cm
\begin{lemma}[see{ [2, Chapter 2, Proposition 2.1]}] Let \,
$\mathcal{G}$\, be an MD-algebra with F $\in {\mathcal{G}}^{*}$ is
not vanishing perfectly in ${\mathcal{G}}^{1}$, i.e. there exists $U
\in {\mathcal{G}}^{1}$ such that $\langle F, U \rangle \neq 0.$ Then
the K-orbit ${\Omega}_{F}$ is  one of the K-orbits having maximal
dimension.
\end{lemma}

\begin{proof} Assume that
${\Omega}_{F}$ is not a K-orbit with maximal dimension, that is,
$dim{\Omega}_{F} = 0$.  Then we have
$$dim{\mathcal{G}}_{F} = dim{\mathcal{G}} - dim{\Omega}_{F} =
dim{\mathcal{G}}.$$ Consequently, $KerB_{F} = {\mathcal{G}}_{F} =
\mathcal{G} \supset {\mathcal{G}}^{1}$ and F is perfectly vanishing
in ${\mathcal{G}}^{1}$. This contradicts the hypotheses of the
Lemma. Therefore, ${\Omega}_{F}$ must be  a K-orbit with maximal
dimension.\end{proof}

\vskip0,8cm
\begin{lemma} Let F be an arbitrary element of \,${\mathcal{G}^{*}}$.
Then\,  $dim{\Omega}_{F} = rank(B)$, where $B = (b_{ij})_{5}: =
(\langle F,[X_{j}, X_{i}]\rangle), 1\leq i, j \leq 5$, is the matrix
of the skew-symmetric bilinear form $B_{F}$ in the basis $(\nobreak
X_{1}, X_{2}, X_{3}, X_{4}, X_{5} \nobreak)$ of \,
$\mathcal{G}$.\end{lemma}

\begin{proof}
Let $U = aX_{1}+bX_{2}+cX_{3}+dX_{4}+eX_{5}\in \mathcal{G}$. Then we
have
\begin{align} \notag {\mathcal{G}}_{F}& = KerB_{F}\\
\notag & = \{U\in \mathcal{G}/ \langle F, [U, X_{i}] \rangle = 0;\,
i = 1, 2, ,3, 4, 5\}.\end{align} By simple computation, we obtain
$$U\in {\mathcal{G}}_{F}\, \Leftrightarrow \, B
\begin{pmatrix} a\\b\\c\\d\\f\end{pmatrix} =
\begin{pmatrix}0\\0\\0\\0\\0\end{pmatrix}.$$
\vskip5mm Hence, $dim{\Omega}_{F} = dim{\mathcal{G}} -
dim{\mathcal{G}}_{F} = rank(B)$.\end{proof}

\vskip0.8cm
\begin{lemma}If $\mathcal{G}$ is a real solvable Lie algebra of dimension 5 with
the first derived ideal ${\mathcal{G}}^{1} \cong {\mathbb{R}}^{4}$
then $\mathcal{G}$ is a MD5-algebra.\end{lemma}

\begin{proof}Let $\mathcal{G}$ be a real solvable Lie algebra with dimension 5
such that ${\mathcal{G}}^{1}$ is the commutative Lie algebra with
dimension 4. Without loss of generality, we may assume that
${\mathcal{G}}^{1} = \mathbb{R}.X_{2} \oplus \mathbb{R}.X_{3} \oplus
\mathbb{R}.X_{4} \oplus \mathbb{R}.X_{5} \equiv {\mathbb{R}}^{4}$,\,
$ad_{X_{1}}= {(a_{ij})}_{4} \in End({\mathcal{G}}^{1}) \equiv
Mat_{4}(\mathbb{R});\, a_{ij} \in \mathbb{R}, 1\leq i, j\leq 4$.

Let $F = \alpha{X_{1}}^{*} + \beta{X_{2}}^{*} + \gamma{X_{3}}^{*} +
\delta{X_{4}}^{*} + \sigma{X_{5}}^{*}\equiv (\alpha, \beta, \gamma,
\delta, \sigma)$ be an arbitrary element from
${\mathcal{G}}^{*}\equiv {\mathbb{R}}^{5}; \alpha, \beta, \gamma,
\delta, \sigma \in \mathbb{R}$. Then, by simple computation, we can
see that the matrix B of the bilinear form $B_{F}$ in the basis
$(\nobreak X_{1}, X_{2}, X_{3}, X_{4}, X_{5} \nobreak)$\,of
$\mathcal{G}$ is a matrix of the following
$$\begin{pmatrix}
0& -\sum_{i=2}^{5} a_{i2}{\alpha}_{i}& -\sum_{i=2}^{5} a_{i3}{\alpha}_{i}
&-\sum_{i=2}^{5} a_{i4}{\alpha}_{i}& -\sum_{i=2}^{5} a_{i5}{\alpha}_{i}\\
\sum_{i=2}^{5} a_{i2}{\alpha}_{i}&0&0&0&0\\
\sum_{i=2}^{5} a_{i3}{\alpha}_{i}&0&0&0&0\\
\sum_{i=2}^{5} a_{i4}{\alpha}_{i}&0&0&0&0\\
\sum_{i=2}^{5} a_{i2}{\alpha}_{i}&0&0&0&0\end{pmatrix}.$$ \vskip5mm
It is now clear that $rank(B)\in \{ 0, 2\}$. Hence, according to
Lemma 2.4, ${\Omega}_{F}$ is the orbit with dimension 0 or 2, i.e.
$\mathcal{G}$ is an MD5-algebra.
\end{proof}

We now  prove the main theorem.\\

\begin{prove} \end{prove}
It is clear that assertion I of Theorem 2.1 holds obviously. We only
need to prove assertion II. Assume that $\mathcal{G}$ is an
indecomposable MD5-algebra with basis $(X_{1}, X_{2}, X_{3}, X_{4},
X_{5} )$ and its first derived ideal $\mathcal{G}^{1}$ is
commutative. Then $dim\, \mathcal{G}^{1} \in \{ 1, 2, 3, 4\}$. In
[16, Theorem 2.1] and [19, Theorem 3.2], the cases had been
considered when $dim \,\mathcal{G}^{1} \in \{ 3, 4\}$. Therefore, we
only need to consider the remaining cases when $dim
\,\mathcal{G}^{1} \in \{ 1, 2\}$. However, for the sake of completeness,
we now consider here all cases.\\

\begin{itemize}
    \item[1.] $dim {\mathcal{G}}^{1} = 1$. Without loss of
    generality , we may assume that ${\mathcal{G}}^{1} = \mathbb{R}.X_{5}\equiv \mathbb{R}$.

            \item[1.1.] Assume that there exists $i \in \{1, 2, 3, 4\}$ with $[X_{i}, X_{5}]\neq 0.$
         Renumber the given basis, if necessary, and we suppose
         that $[X_{4}, X_{5}]= aX_{5}$, for some $a \in \mathbb{R}\backslash\{0\}$.
         Then, by changing $X_{4}$ with ${X_{4}}^{'} =
            \frac{1}{a}X_{4}, $ we obtain  $[{X_{4}}^{'}, X_{5}]= X_{5}$.
         Now, without any restriction of generality, we can assume
         that $[X_{4}, X_{5}]= X_{5}$.

            Let $[X_{i}, X_{5}]= a_{i}X_{5}, [X_{i}, X_{4}]=
            b_{i}X_{5}; a_{i}, b_{i} \in \mathbb{R}; i = 1, 2, 3$.
            Then, by changing ${X_{i}}^{'} = X_{i} - a_{i}X_{4} + b_{i}X_{5} (i = 1, 2, 3)$,
            we get $[{X_{i}}^{'}, X_{5}]= [{X_{i}}^{'}, X_{4}]= 0; i = 1, 2, 3.$.
            Hence, we can always suppose right from the start that
            $[X_{i}, X_{5}]= [X_{i}, X_{4}]= 0; i = 1, 2, 3$.

            Now, let $[X_{i}, X_{j}]= c_{ij}X_{5}, \, c_{ij}\in \mathbb{R};
            1 \leq i < j \leq 3$. Then, by using the Jacobi identity,  we
            get $c_{ij} = 0$ for all $i, j,\, 1 \leq i < j \leq 3$.
            But this shows that $\mathcal{G}$ is decomposable, which is  a
            contradiction. Thus, this case cannot happen.

            \item[1.2.] Assume that $[X_{i}, X_{5}]= 0$ for all i = 1, 2, 3, 4.
            Then, there exists $[X_{i}, X_{j}]= c_{ij}X_{5},\, c_{ij}\neq
            0$ for some $i, j  \in \{1, 2, 3, 4\}, i \neq j$.
            By applying the same argument as in Case 1.1, we can suppose that
            $[X_{1}, X_{2}]= [X_{3}, X_{4}] = X_{5}$ and $[X_{i}, X_{3}]= [X_{i},
            X_{4}]= 0; i = 1,2$. Therefore, $\mathcal{G}\cong
            \mathcal{G}_{5,1}$.\\

          \item[2.] $dim {\mathcal{G}}^{1} = 2$. Without
           loss of generality, we now assume that
           ${\mathcal{G}}^{1} = \mathbb{R}.X_{4} \oplus \mathbb{R}.X_{5}\equiv {\mathbb{R}}^{2}$;\,
           $ad_{X_{1}},\, ad_{X_{2}},\, ad_{X_{3}} \in End({\mathcal{G}}^{1})
           \equiv Mat_{2}(\mathbb{R})$.

            \item[2.1.]$[X_{i}, X_{j}]= 0, 1\leq i < j \leq 3.$

            If there exists $ad_{{X}_{i}} = 0$ then $\mathcal{G}$ is
            decomposable, which is a contradiction. Hence, $ad_{{X}_{i}} \neq
            0, i =1, 2, 3$. We now show that we can always
            obtain  $ad_{{X}_{2}} = 0$ by changing the basis. Indeed, we can let $ad_{{X}_{i}}$
            be $\begin{pmatrix} a_{i}&b_{i}\\c_{i}&d_{i}\end{pmatrix}
            \neq \begin{pmatrix} 0&0\\0&0\end{pmatrix}; \, \, i = 1, 2, 3$.
            We first assume that $a_{3} \neq 0$. Then, by writing
            ${X_{i}}^{'} = X_{i} - \frac{a_{i}}{a_{3}}X_{3} $,
            we get $ad_{{X}_{i}^{'}} =
            \begin{pmatrix} 0&{b_{i}}^{'}\\{c_{i}}^{'}&{d_{i}}^{'}\end{pmatrix},
            \, i = 1, 2$. Hence, we can suppose  that
            $ad_{{X}_{i}} =
            \begin{pmatrix} 0&b_{i}\\c_{i}&d_{i}\end{pmatrix}, \, i = 1,
            2$. According to Lemma 2.3, $ad_{X_{1}}\circ ad_{X_{2}}
            = ad_{X_{2}}\circ ad_{X_{1}}$. It follows that
            $ad_{{X}_{1}} = k.ad_{{X}_{2}},$ for some $k \in
            \mathbb{R}\backslash \{0\}$. By changing ${X_{2}}^{'} = X_{2} - k.X_{1} $,
            we get $ad_{{X}_{2}^{'}} = 0$, a contradiction. When $d_{3} \neq
            0$, then by using the same argument, we can also obtain a
            contradiction. Finally, assume that $a_{3} = d_{3} = 0$,
            $b_{3}^{2} + c_{3}^{2} \neq 0$. In view of Lemma 2.3, we
            get $ad_{X_{i}}\circ ad_{X_{3}} = ad_{X_{3}}\circ ad_{X_{i}}(i = 1, 2)$

            Hence, it  follows that $ad_{{X}_{i}} = k_{i} \begin{pmatrix}1&0\\0&1
            \end{pmatrix}$  $0\neq k_{i} \in \mathbb{R},\, i = 1, 2$. In particular,
            $ad_{{X}_{2}} = k.ad_{{X}_{1}}, k = \frac{k_{2}}{k_{1}}$.
            Now, by changing ${X_{2}}^{'} = X_{2} - k.X_{1} $,
            we get $ad_{{X}_{2}^{'}} = 0$, again a contradiction.
            Hence, Case 2.1 can not happen.

           \item[2.2.] Assume that there exists $[X_{i}, X_{j}]\neq 0, 1\leq i < j \leq
            3$ and $ad_{{X}_{i}} = 0$ \, ${i =1, 2, 3}$.

            It is clear that $\mathcal{G}^{1} = \langle [X_{1}, X_{2}], [X_{1}, X_{3}], [X_{2},
            X_{3}]\rangle$, and whence, the rank of $\{[X_{1}, X_{2}], [X_{1}, X_{3}], [X_{2},
            X_{3}]\}$ is 2 and without restriction of generality, we
            can assume that $\{[X_{1}, X_{2}], [X_{2}, X_{3}]\}$ is
            a basis of $\mathcal{G}^{1}$. Let $[X_{1}, X_{2}] = aX_{4} +
            bX_{5}$, $[X_{2}, X_{3}] = cX_{4} + dX_{5}$
            with $D: = det\begin{pmatrix}a&b\\c&d\end{pmatrix}\neq 0$.
            By changing basis as follows
            $$X_{4} = \frac{1}{D} (d{X_{4}}^{'} - b{X_{5}}^{'}),\,
            X_{5} = \frac{1}{D} (- c{X_{4}}^{'} + a{X_{5}}^{'})$$
            we get $[X_{1}, X_{2}] = {X_{4}}^{'},\, [X_{2}, X_{3}] =
            {X_{5}}^{'}$. Hence, we can assume
            that $[X_{1}, X_{2}] = X_{4},\, [X_{2}, X_{3}] = X_{5}$.

            Let $[X_{1}, X_{3}] = \alpha X_{4} + \beta X_{5}$.
            Then, by changing the  basis as follows:
            $${X_{1}}^{'} = X_{1} - \beta X_{2}, {X_{2}}^{'} = X_{2}, {X_{3}}^{'} = -\alpha X_{2} + X_{3}$$
            we get
            $$[{X_{1}}^{'}, {X_{2}}^{'}] = X_{4},\, [{X_{2}}^{'}, {X_{3}}^{'}] = X_{5},
            \,[{X_{1}}^{'}, {X_{3}}^{'}] = 0. $$
            Thus, we can always assume that
            $$[X_{1}, X_{2}] = X_{4},\, [X_{2}, X_{3}] = X_{5},\,[X_{1}, X_{3}] = 0. $$
            Therefore \, $\mathcal{G}\cong {\mathcal{G}}_{5,2,1}$.

            \item[2.3.]Assume that there exists $[X_{i}, X_{j}]\neq 0$ and $ad_{{X}_{k}} \neq
            0$, $ 1\leq i \neq j \leq 3,\, 1\leq k \leq 3 $. Then, without
            loss of generality, we may assume that $ad_{{X}_{3}} \neq 0$.

            We can always change basis of $\mathcal{G}^{1}$ such that $ad_{{X}_{3}}$
            becomes one of the following matrices
            $$\begin{pmatrix}0&0\\1&0\end{pmatrix},
            \begin{pmatrix}1&0\\0&\lambda\end{pmatrix},
            \begin{pmatrix}1&1\\0&1\end{pmatrix},
            \begin{pmatrix}cos\varphi& -sin\varphi \\sin\varphi&cos\varphi
            \end{pmatrix};\, \lambda \in \mathbb{R}, \varphi \in (0, \pi).$$

             \item[2.3a.]Assume that $ad_{X_{3}} =
                \begin{pmatrix}0&0\\1&0\end{pmatrix}$.
                Then by using an argument analogous to that in Subsection 2.2, we get
                $ad_{X_{1}} = ad_{X_{2}} = 0$. Again, by  Jacobi identity, we obtain
                $[X_{1}, X_{2}] = a X_{5},\, a \in \mathbb{R}$.

                Let $[X_{i}, X_{3}] = a_{i}X_{4} + b_{i} X_{5};\, a_{i}, b_{i} \in \mathbb{R},\, i = 1, 2$.
                If $a = 0$, then by changing ${X_{i}}^{'} = X_{i} + b_{i}X_{4},$ we get
                $[{X_{i}}^{'}, X_{3}] = a_{i}X_{4},\, i =1, 2$.
                Hence, we can always assume from the outset that $[X_{i}, X_{3}] =
                a_{i}X_{4};\, i = 1, 2;\, {a_{1}}^{2} + {a_{2}}^{2} \neq 0$.
                Without loss of generality, we may assume
                that $a_{2} \neq 0$. Now, we change again the basis as follows
                $${X_{1}}^{'} = X_{1} - \frac{a_{1}}{a_{2}}X_{2},\, {X_{2}}^{'} =
                \frac{1}{a_{2}}X_{2}.$$
                Then we get $[{X_{1}}^{'}, X_{3}] = 0, [{X_{2}}^{'}, X_{3}] =
                X_{4}$, i.e. $\mathcal{G}$ is decomposable, a
                contradiction. Hence, $a \neq 0$.\\

                In the same way, we obtain
                $$[X_{1}, X_{2}] = [X_{3}, X_{4}] = X_{5},\,[X_{2}, X_{3}] = \lambda X_{4},
                \,0 \neq \lambda \in \mathbb{R}.$$
                Therefore\, $\mathcal{G} \cong {\mathcal{G}}_{5,2,2(\lambda)}$.\\
                \item[2.3b.]Assume that $ad_{X_{3}} =
                \begin{pmatrix}1&0\\0&\lambda\end{pmatrix},\, \lambda \in \mathbb{R}$.
                Then, by using a similar argument as above, we get ${\mathcal{G}}_{5,2,3}$:
                $[X_{1}, X_{2}] = X_{5}, [X_{3}, X_{4}] = X_{4}$.
                By using Lemmas 2.2\, and\, 2.3, and by direct computation, we can  show that
                ${\mathcal{G}}_{5,2,3}$ is not an MD5-algebra. Hence, this case has to be rejected.\\

                \item[2.3c.] Assume that $ad_{X_{3}} \in \left\{ \begin{pmatrix}1&1\\0&1\end{pmatrix},
                \begin{pmatrix}cos\varphi& -sin\varphi \\sin\varphi&cos\varphi
                \end{pmatrix};\,\varphi \in (0,\pi)\right\}$. By using a similar  argument as above,
                these cases have to be also rejected.\\

              \item[3.] $dim {\mathcal{G}}^{1} = 3$. We can
              always change basis to obtain ${\mathcal{G}}^{1} =
             \mathbb{R}.X_{3} \oplus \mathbb{R}.X_{4} \oplus
             \mathbb{R}.X_{5} \equiv {\mathbb{R}}^{3}$;\,
             $ad_{X_{1}},\, ad_{X_{2}} \in End({\mathcal{G}}^{1})
            \equiv Mat_{3}(\mathbb{R})$.\\
\end{itemize}

It is obvious that $ad_{X_{1}}$ and $ad_{X_{2}}$ cannot be
 the trivial operators concurrently because ${\mathcal{G}}^{1} \cong
{\mathbb{R}}^{3}$. Without loss of generality, we may assume that
$ad_{X_{2}} \neq 0$. Then, by changing basis, if necessary, we obtain a
similar classification of $ad_{{X}_{2}}$ as follows\\
\begin{itemize}
\item $\begin{pmatrix} {{\lambda}_1}&0&0\\
0&{{\lambda}_2}&0\\0&0&1 \end{pmatrix},\quad({\lambda}_1,
{\lambda}_2 \in \mathbb{R}\setminus \lbrace 1\rbrace, \, {\lambda}_1
\neq {\lambda}_2 \neq 0$);
\item $\begin{pmatrix}1&0&0\\0&1&0\\0&0&{\lambda}
\end{pmatrix},\quad ({\lambda} \in \mathbb{R}\setminus \lbrace 0, 1
\rbrace$);
\item $\begin{pmatrix} {\lambda}&0&0\\0&1&0\\0&0&1
\end{pmatrix},\quad ({\lambda} \in \mathbb{R}\setminus \lbrace 1
\rbrace$);
\item $\begin{pmatrix} 1&0&0\\0&1&0\\0&0&1
\end{pmatrix}$;
\item $\begin{pmatrix} {\lambda}&0&0\\0&1&1\\0&0&1
\end{pmatrix},\quad ({\lambda} \in \mathbb{R}\setminus \lbrace 1
\rbrace$);
\item $\begin{pmatrix} 1&1&0\\0&1&0\\0&0&{\lambda}
\end{pmatrix},\quad ({\lambda} \in \mathbb{R}\setminus \lbrace 0, 1
\rbrace$);
\item $\begin{pmatrix} 1&1&0\\0&1&1\\0&0&1
\end{pmatrix}$;
\item $\begin{pmatrix} cos{\varphi}&-sin{\varphi}&0\\
sin{\varphi}&cos{\varphi}&0\\0&0&\lambda \end{pmatrix},\quad
(\lambda \in \mathbb{R}\setminus \lbrace 0\rbrace, \, \varphi \in
(0, \pi))$.\\
\end{itemize}

Assume that $[X_{1}, X_{2}] = mX_{3} + nX_{4} + pX_{5}; m, n, p \in
\mathbb{R}$. We can always change basis to have $[X_{1}, X_{2}] =
mX_{3}$. Indeed, if
 $$ad_{X_{2}} = \begin{pmatrix} {{\lambda}_1}&0&0\\
0&{{\lambda}_2}&0\\0&0&1 \end{pmatrix},\, ({\lambda}_1, {\lambda}_2
\in \mathbb{R}\setminus \lbrace 1\rbrace, \, {\lambda}_1 \neq
{\lambda}_2 \neq 0),$$ then by changing $X_{1}$ for ${X_{1}}^{'} =
X_{1} + \frac{n}{\lambda_{2}}X_{4} + pX_{5}$ we get $[{X_{1}}^{'},
X_{2}] = mX_{3}$, $m \in \mathbb{R}$. For the other values of
$ad_{X_{2}}$, we can also change basis in the same way. Hence,
without restriction of generality, we can assume that $[X_{1},
X_{2}] = mX_{3}$, $m \in \mathbb{R}$.

There are three cases which contradict each other as follows.
\begin{itemize}
    \item[3.1.] $[X_{1}, X_{2}] = 0$ ( i.e. $m = 0$ ) and
    $ad_{X_{1}} = 0$. Then $\mathcal{G} =
\mathcal{H} \oplus \mathbb{R}.X_{1}$, where $\mathcal{H}$ is the
subalgebra of $\mathcal{G}$ generated by $\{X_{2}, X_{3}, X_{4},
X_{5}\}$, i.e. $\mathcal{G}$ is decomposable. Hence, this case is
rejected.
    \item[3.2.] $[X_{1}, X_{2}] = 0$ and $ad_{X_{1}} \neq 0$.\\
    \item[3.2a.] Assume that $ad_{X_{2}} = \begin{pmatrix} {{\lambda}_1}&0&0\\
0&{{\lambda}_2}&0\\0&0&1 \end{pmatrix}; \, \,{\lambda}_1,\,
{\lambda}_2 \in \mathbb{R}\setminus \lbrace 1\rbrace,\, \,
{\lambda}_1 \neq {\lambda}_2 \neq 0.$ In view of Lemma 2.1, it
follows by a direct computation that
$$ad_{X_{1}} = \begin{pmatrix} \mu&0&0\\0&\nu&0 \\ 0&0&\xi
\end{pmatrix};\, \mu, \nu, \xi \in \mathbb{R};\, {\mu}^{2} + {\nu}^{2} +{\xi}^{2} \ne
0.$$ If $\xi \neq 0$, by changing ${X_{1}}^{'} = X_{1} - \xi X_{2}$,
we get $$ad_{{X_{1}}^{'}} =  \begin{pmatrix} {\mu}^{'}&0&0\\0&{\nu}^{'}&0 \\
0&0&0 \end{pmatrix};$$ where \, ${\mu}^{'} = \mu - \xi
{\lambda}_{1}, {\nu}^{'} = \nu - \xi {\lambda}_{2}.$ Thus, we can
assume that $$ad_{X_{1}} = \begin{pmatrix}
\mu&0&0\\0&\nu&0\\ 0&0&0 \end{pmatrix};\, \mu, \nu  \in
\mathbb{R};\, {\mu}^{2} + {\nu}^{2} \ne 0.$$ Using Lemmas 2.2, 2.3,
and by direct computation, we can  show that $\mathcal{G}$ will not be an
MD5-algebra in Case 3.2a . So this case must be rejected.\\
\item[3.2b.] In exactly the same way, but replacing the
considered value of $ad_{X_{2}}$ with the others, we can easily see
that Case 3.2 cannot occur.\\
\item[3.3.] $[X_{1}, X_{2}] \neq 0$ ( i.e. $m \neq 0$ ).
By changing $X_{1}$ by ${X_{1}}^{'} = {\frac{1}{m}}X_{1},$ we have
$[{X_{1}}^{'}, X_{2}] = X_{3}$. Hence, without loss of generality,
we may assume that $[X_{1}, X_{2}] = X_{3}$. By using a similar
argument as the one in Case 3.2a, we obtain again a contradiction if
$ad_{X_{1}} \neq 0$. In other words, $ad_{X_{1}} = 0$. Therefore, in
the dependence on the value of $ad_{X_{2}}$,\,\, $\mathcal{G}$ must
be isomorphic to one of the following algebras:

\begin{itemize}
\item ${\mathcal{G}}_{5,3,1({\lambda}_{1}, {\lambda}_{2})},\quad
({\lambda}_{1}, {\lambda}_{2} \in \mathbb{R}\setminus \{1\},
{\lambda}_{1} \neq {\lambda}_{2} \neq 0 )$;
\item ${\mathcal{G}}_{5,3,2(\lambda)},\quad (\lambda \in \mathbb{R} \setminus
\{0, 1\})$;
\item ${\mathcal{G}}_{5,3,3(\lambda)},\quad(\lambda \in
\mathbb{R}\nobreak \setminus \{1\})$;
\item ${\mathcal{G}}_{5,3,4}$;
\item ${\mathcal{G}}_{5,3,5(\lambda)},\quad (\lambda \in
\mathbb{R}\setminus \{1\})$;
\item ${\mathcal{G}}_{5,3,6(\lambda)},\quad(\lambda \in \mathbb{R}\setminus
\{0, 1\})$;
\item ${\mathcal{G}}_{5,3,7}$;
\item ${\mathcal{G}}_{5,3,8(\lambda, \varphi)},\quad (\lambda \in
\mathbb{R}\setminus \{0\}),  \,  \varphi \in (0, \pi)\,)$.
\end{itemize}
Obviously, these algebras are not mutually isomorphic to each
other.\\
\item[4.] $dim {\mathcal{G}}^{1} = 4$. Without
loss of generality, we may assume that ${\mathcal{G}}^{1} =
\mathbb{R}.X_{2} \oplus \mathbb{R}.X_{3} \oplus \mathbb{R}.X_{4}
\oplus \mathbb{R}.X_{5} \equiv {\mathbb{R}}^{4}$,\, $ad_{X_{1}} \in
End({\mathcal{G}}^{1}) \equiv Mat_{4}(\mathbb{R})$.
\end{itemize}

According to Lemma 2.5, the final assertions of Theorem 2.1 can be
obtained by using similar classification of $ad_{{X}_{1}}$.
\vskip0.5cm
In view of Lemma 2.4, it follows by  direct computation
that all algebras listed in Theorem 2.1 are MD5-algebras. This
completes the
proof. \hfill $\square$\\

\subsection*{Concluding Remark} Recall that every real Lie
algebra $\mathcal{G}$ defines only one connected and simply
connected Lie group G such that Lie(G) = $\mathcal{G}$. Therefore,
we obtain a collection of twenty - five families of connected and
simply connected MD5-groups corresponding to given indecomposable
MD5-algebras in Theorem 2.1. For the sake of convenience, we denote
every MD5-group from this collection by using the same indices as
its corresponding MD5-algebra. For example, $G_{5,3,1({\lambda}_{1},
{\lambda}_{2})}$ is the connected and simply connected MD5-group
which corresponds to ${\mathcal{G}}_{5,3,1({\lambda}_{1},
{\lambda}_{2})}$. All of these groups are indecomposable MD5-groups.
In the next papers, we shall compute the invariants of given
MD5-algebras, describe the geometry of K-orbits of its corresponding
MD5-groups and also we shall classify topologically the
MD5-foliations associated with these MD5-groups. In addition,
characterization theorems of Connes $C^{*}$-algebras corresponding
to these MD5-foliations will also be established.\\

\subsection*{Acknowledgement}
The first author would like to thank Professor Do Ngoc Diep for
giving him excellent advice and support. He also want to thank the
Organizing Committee of The Second International Congress In
Algebras and Combinatorics - July 2007, Beijing, China for inviting
him to give a talk on this topic at the congress.\\

\end{document}